\documentclass[12pt]{article}
\usepackage{graphicx,amssymb,amsmath,vmargin,multicol}
\usepackage[english]{babel}
\usepackage[T1]{fontenc}
\usepackage[latin1]{inputenc}
\catcode`\«=\active
\catcode`\»=\active
\def«{\og\ignorespaces}%
\def»{{\fg}}%
\title{$L^p$ version of a result by Rankin\cite{rankin}}
\setpapersize{custom}{21cm}{29.7cm}
\setmarginsrb{1cm}{1cm}{1cm}{2cm}{0pt}{0pt}{0pt}{0pt}

\author{Mathieu Dutour \footnote{Mathieu.Dutour@wanadoo.fr, 26 rue de la république, 45000 Orléans France}}
\date{Jul 2000}
\begin{document}
\def\Tore#1{\leavevmode\kern-.0em\raise.2ex\hbox{$\R^2$}\kern-.1em/\kern-.1em\lower.25ex\hbox{$#1$}}

\newcommand{\D}{\displaystyle}
\newcommand{\R}{\ensuremath{\Bbb{R}}}
\newcommand{\N}{\ensuremath{\Bbb{N}}}
\newcommand{\Q}{\ensuremath{\Bbb{Q}}}
\newcommand{\C}{\ensuremath{\Bbb{C}}}
\newcommand{\Z}{\ensuremath{\Bbb{Z}}}
\newcommand{\T}{\ensuremath{\Bbb{T}}}
\newcommand{\K}{\ensuremath{\Bbb{K}}}
\newcommand{\curl}{\ensuremath{\mathop{\rm curl}\nolimits}}
\newcommand{\divergence}{\ensuremath{\mathop{\rm div}\nolimits}}
\newcommand{\determinant}{\ensuremath{\mathop{\rm det}\nolimits}}
\newcommand{\Rez}{\ensuremath{\mathop{\rm Re}\nolimits}}
\newcommand{\Imz}{\ensuremath{\mathop{\rm Im}\nolimits}}

\newtheorem{theorem}{Theorem}[section]
\newtheorem{lemme}[theorem]{Lemma}
\newtheorem{proposition}[theorem]{Proposition}
\newtheorem{definition}[theorem]{Definition}
\newtheorem{remarque}[theorem]{Remarque}
\newtheorem{corollaire}[theorem]{Corollaire}
\newtheorem{exemple}[theorem]{Exemple}
\newtheorem{probleme}[theorem]{Problème}
\newtheorem{hypoth}[theorem]{Hypothèse}
\maketitle
\begin{abstract}
\noindent We extend a classical result by Rankin\cite{rankin}. We consider the following question: given $n$ vectors $v_i$ in the ball of radius $R$ of an infinite dimensional Banach space ${\cal B}$ with $d(v_i,v_j)\geq 1$, can we bound the number $n$?
\end{abstract}

\section{Packing Problem}
\noindent Here is a classical result by Rankin, cited in \cite{deza} p.~85
\begin{theorem}
Let $R<\frac{1}{\sqrt{2}}$ and let $N_R$ denote the maximum number $N$ of points $x_1, x_2,\dots, x_N$ that can be placed in a closed (euclidean) ball of radius $R$ in such a way that $||x_i-x_j||\geq 1$ for all $i\not=j=1,\dots,N$ then $N_R=\lfloor \frac{1}{1-2R^2}\rfloor$.
\end{theorem}
\noindent Our goal is to extend this result to infinite dimensional Banach space ${\cal B}$. We will solve the problem in $l_p$.\\
We denote by $B_R$ the ball of center $0$ and radius $R$ for the $||.||$ norm of ${\cal B}$.\\
We define 
\begin{equation}
{\cal E}_{{\cal B},n,R}=
\left\lbrace\begin{array}{c}
(v_1,v_2,\dots,v_n)\in B_R,\\
||v_i-v_j||\geq 1\mbox{~if~}i\not=j
\end{array}\right\rbrace
\end{equation}
If $R\geq \frac{1}{2}$ then ${\cal E}_{{\cal B},n,R}\not=\emptyset$ for all $n\in\N$.\\
\begin{equation}\label{trivial-bound}
{\cal E}_{{\cal B},2,R}=\emptyset\Leftrightarrow R<\frac{1}{2}
\end{equation}
We define also
\begin{equation}
N({\cal B},R)=\sup\,\{n\in\N,\,\,\mbox{~such~that~}{\cal E}_{{\cal B},n,R}\not=\emptyset\}
\end{equation}
and
\begin{equation}
R_c({\cal B})=\inf\{R\in\R_+\mbox{~such~that~} N({\cal B},R)=\infty\}
\end{equation}
We will have then
\begin{enumerate}
\item If $R<R_c({\cal B})$ then $n\leq N({\cal B},R)<\infty$ if $(v_1,v_2,\dots,v_n)\in {\cal E}_{{\cal B},n,R}$.
\item If $R>R_c({\cal B})$ then there is no way to majorize the number $n$.
\end{enumerate}
From \ref{trivial-bound} we know that $R_c({\cal B})\geq \frac{1}{2}$ for any infinite dimensional Banach space ${\cal B}$.

\section{A linear mapping}
\noindent Let us denote by $E=l_p$ the Banach subspace of $\R^{\N}$ equipped with the norm
\begin{equation}
\begin{array}{rcl}
N_p:E&\mapsto&\R\\
x&\mapsto&\sqrt[p]{\sum_{i=1}^{\infty}|x_i|^p}
\end{array}
\end{equation}
We consider the linear mapping
\begin{equation}
\begin{array}{rcl}
\phi:E^n&\mapsto&E^{\frac{n(n-1)}{2}}\\
(v_1,v_2,\dots,v_n)&\mapsto&(v_i-v_j)_{1\leq i<j\leq n}
\end{array}
\end{equation}
We equip $E^s$ with the norm
\begin{equation}
\begin{array}{rcl}
N_{p,s}:E^s&\mapsto&\R\\
(v_1,v_2,\dots,v_s)&\mapsto&\sqrt[p]{\sum_{i=1}^{s}N_p^p(v_i)}
\end{array}
\end{equation}
Question is to evaluate the norm of the mapping $\phi$ with $E^n$ and $E^{\frac{n(n-1)}{2}}$ equipped with the norms $N_{p,n}$ and $N_{p,\frac{n(n-1)}{2}}$.\\
We denote by $||\phi||_{p}$ the norm of this mapping
\begin{theorem}
We have $||\phi||_1=n-1$
\end{theorem}
{\it Proof} What we need to find is the best constant $C$ such that
\begin{equation}
\sum_{1\leq i<j\leq n}N_1(v_i-v_j)\leq C\sum_{i=1}^{n} N_1(v_i)
\end{equation}
Let us majorize the sum
\begin{equation}
\begin{array}{rcl}
\sum_{1\leq i<j\leq n}N_1(v_i-v_j)
&\leq&\sum_{1\leq i<j\leq n}N_1(v_i)+N_1(v_j)\\
&\leq &(n-1)\sum_{i=1}^{n}N_1(v_i)
\end{array}
\end{equation}
So we have constant $C\leq n-1$. Then if we set
\begin{equation}
(v_1,v_2,\dots,v_n)=(h,0,\dots,0)
\end{equation}
we obtain equality.

\begin{theorem}
We have $||\phi||_2=\sqrt{n}$
\end{theorem}
{\it Proof} A well known formula (Maybe RANKIN, but needs to be verified) is
\begin{equation}
\sum_{1\leq i<j\leq n}N_2(x_i-x_j)^2=n\sum_{i=1}^{n}N_2(x_i)^2-N_2(\sum_{i=1}^{n}x_i)^2 ,
\end{equation}
from which we get 
\begin{equation}
\sqrt{\sum_{1\leq i<j\leq n}N_2(x_i-x_j)^2}\leq \sqrt{n}\sqrt{\sum_{i=1}^{n}N_2(x_i)^2} ,
\end{equation}
which provides $||\phi||_2\leq \sqrt{n}$. Equality is attained with
\begin{equation}
v=(v_1,\dots,v_n)\mbox{~with~}\sum v_i=0
\end{equation}

\begin{theorem}
We have $||\phi||_{\infty}=2$
\end{theorem}
{\it Proof} What we need to find is the best constant $C$ such that
\begin{equation}
\sup_{1\leq i<j\leq n}N_{\infty}(v_i-v_j)\leq C\sup_{1\leq i\leq n} N_{\infty}(v_i)
\end{equation}
Let us majorize the sup
\begin{equation}
\begin{array}{rcl}
\sup_{1\leq i<j\leq n}N_{\infty}(v_i-v_j)
&\leq&\sup_{1\leq i<j\leq n}N_{\infty}(v_i)+N_{\infty}(v_j)\\
&\leq &2\sup_{1\leq i\leq n} N_{\infty}(v_i)
\end{array}
\end{equation}
which gives $||\phi||_{\infty}\leq 2$. Equality is obtained with
\begin{equation}
v=(h,-h,0,\dots,0)
\end{equation}

We use a classic interpolation theorem
\begin{theorem}\label{interpolation-theorem}
({\it Stein Interpolation theorem (see \cite{reedsimonII}, p.~40)}) If $p\leq q\leq r$ and $\phi:L^h\mapsto L^h$ is a continous mapping for $h=p,r$ then $\phi$ is continuous for $q$ and
\begin{equation}
||\phi||_{L^q} < [||\phi||_{L^r}]^{\frac{(q-p)r}{(r-p)q}}   [||\phi||_{L^p}]^{\frac{(r-q)p}{(r-p)q}}
\end{equation}
\end{theorem}
this theorems gives, if $1<q<2$, 
\begin{equation}
\begin{array}{rcl}
||\phi||^q_{q}
&\leq&[||\phi||_2]^{(q-1)2}     [||\phi||_1]^{(2-q)}\\
&\leq&[\sqrt{n}]^{(q-1)2}     [n-1]^{(2-q)}\\
&\leq&[n]^{q-1}     [n-1]^{2-q}
\end{array}
\end{equation}
It gives also, if $1<q<\infty$,
\begin{equation}
\begin{array}{rcl}
||\phi||^q_{q}
&\leq&[||\phi||_{\infty}]^{q-2}     [||\phi||_2]^{2}\\
&\leq&[2]^{q-2}     [\sqrt{n}]^{2}\\
&\leq&n[2]^{q-2}
\end{array}
\end{equation}

\section{Lower Bound}\label{sec:lower-bound}
\begin{theorem}
If $1<q<2$, $R<2^{-\frac{1}{q}}$ and ${\cal E}_{l_q,n,R}\not= \emptyset$ then
\begin{equation}\label{bound-ONE}
n\leq \lfloor\,\,\frac{1}{1-[2R^q]^{\frac{1}{q-1}}}\,\,\rfloor=\psi(q,R)
\end{equation}
\end{theorem}
{\it Proof} Assume $(v_1,\dots,v_n)$ belongs to ${\cal E}_{l_q,n,R}$ we then obtain by the linear application bound
\begin{equation}
\sum_{1\leq i<j\leq n}N_q(v_i-v_j)^q\leq u(q)^q\sum_{i=1}^{n}N_q(v_i)^q
\end{equation}
and so
\begin{equation}
\frac{n(n-1)}{2}\leq u(q)^qnR^q
\end{equation}
and this gives us
\begin{equation}
\frac{1}{2R^q}\leq (\frac{n}{n-1})^{q-1}
\end{equation}
which inverts into
\begin{equation}
n\leq \frac{1}{1-[2R^q]^{\frac{1}{q-1}}}
\end{equation}
and we have the result.

\begin{theorem}
If $1<q$, $R<2^{\frac{1}{q}-1}$ and ${\cal E}_{l_q,n,R}\not=\emptyset$ then
\begin{equation}\label{bound-TWO}
n\leq \lfloor\,\,\frac{1}{1-2^{q-1}R^q}\,\,\rfloor=\psi(q,R)
\end{equation}
\end{theorem}
{\it Proof} Assume $(v_1,\dots,v_n)$ belongs to ${\cal E}_{l_q,n,R}$ we then obtain by the linear application bound
\begin{equation}
\sum_{1\leq i<j\leq n}N_q(v_i-v_j)^q\leq u(q)^q\sum_{i=1}^{n}N_q(v_i)^q
\end{equation}
and so
\begin{equation}
\begin{array}{rcl}
\frac{n(n-1)}{2}
&\leq &nR^q ||\phi||_q^q\\
&\leq &n[2]^{q-2}   nR^q\\
&\leq &n^2[2]^{q-2} R^q
\end{array}
\end{equation}
We then obtain
\begin{equation}
n-1\leq n 2^{q-1} R^q
\end{equation}
and so
\begin{equation}
n(1- 2^{q-1} R^q)\leq 1
\end{equation}
which gives us
\begin{equation}
n\leq \lfloor\,\,\frac{1}{1-2^{q-1}R^q}\,\,\rfloor
\end{equation}

\section{Upper bound}\label{sec:upper-bound}
\begin{theorem}
We have $R_c(l_p)\leq 2^{-\frac{1}{p}}$.
\end{theorem}
{\it Proof} It suffices to find vectors having pairwise $L^p$-distance equal to $1$ and norm equal to $2^{-\frac{1}{p}}$.\\
We define
\begin{equation}
e_{i}=\frac{1}{\sqrt[p]{2}}(0,\dots,0,1,0,\dots,0,\dots)\in l_p \mbox{~with~}i\in\N
\end{equation}
and we have
\begin{equation}
N_p(e_{i})=\frac{1}{\sqrt[p]{2}}\mbox{~and~}N_p(e_i-e_j)=1,\mbox{~if~}i\not= j
\end{equation}
since the number of vectors $e_i$ is arbitrary we have the result.

\begin{theorem}
We have $R_c(l_p)\leq 2^{\frac{1}{p}-1}$.
\end{theorem}
{\it Proof} It suffices to find vectors having pairwise $L^p$-distance equal to $1$ and norm equal to $2^{\frac{1}{p}-1}$.\\
Let $n$ such that a Hadamard matrix $A=(a_{ij})\in M_n(\R)$ exists. We write $n=2h$ and define vectors
\begin{equation}
f_{i}=\frac{1}{\sqrt[p]{h}}\frac{1}{2}(a_{i1},a_{i2},\dots,a_{in},0,\dots,0,\dots)\in l_p\mbox{~with~}1\leq i\leq n
\end{equation}
we have
\begin{equation}
N_p(f_{i})=\frac{1}{2}2^{\frac{1}{p}}
\end{equation}
if $i\not= j$ then by orthogonality the number of different coefficient between $f_i$ and $f_j$ is $h$ and we then obtain
\begin{equation}
N_p(f_{i}-f_{j})=1
\end{equation}
If $n=2^r$ then there exist an Hadamard matrix. So the number of vectors having pairwise $L^p$-distance equal to $1$ and norm equal to $2^{\frac{1}{p}-1}$ is not bounded and we get the result: $R_c(l_p)\leq 2^{\frac{1}{p}-1}$.\\
\\
Combining the preceding results we conclude
\begin{equation}
R_c(l_p)=\left\lbrace\begin{array}{rcl}
2^{-\frac{1}{p}}&\mbox{~if~}&1\leq p\leq 2\\
2^{\frac{1}{p}-1}&\mbox{~if~}&2\leq p\leq \infty
\end{array}\right.
\end{equation}
and we remark that $R_c(l_p)=R_c(l_q)$ if $\frac{1}{p}+\frac{1}{q}=1$ which suggest our result has some link with duality.\\
Even more there is a duality in the bound $\psi$ defined implicitely at \ref{bound-ONE} and \ref{bound-TWO}: $\psi(p,R)=\psi(q,R)$ if $\frac{1}{p}+\frac{1}{q}=1$.\\
This result is not new, it seems to appear on (\cite{thanks-maurey}, p.~31-34).\\
I thank Bernard Maurey for useful comment on this paper.

\section{An uniform bound}
\noindent We use here Dvoretsky result \cite{pisier1} to proove that $R_c({\cal B})\leq \frac{1}{\sqrt{2}}$.
\begin{theorem}
Let ${\cal B}$ an infinite dimensional Banach space with norm $||.||$. $\forall \epsilon>0$, $\forall n\in\N^*$ there exists $x_1,x_2,\dots,x_n\in {\cal B}$ such that
\begin{equation}
\forall \alpha\in\R^n,\,\,(1-\epsilon)\sqrt{\sum_{i=1}^{i=n}\alpha_i^2}\leq ||\sum_{i=1}^{i=n}\alpha_i x_i||\leq (1+\epsilon)\sqrt{\sum_{i=1}^{i=n}\alpha_i^2}
\end{equation}
\end{theorem}
Setting $w_i=\frac{1}{\sqrt{2}(1-\epsilon)}$ we have $||w_i-w_j||\geq 1$ and $||w_i||\leq \frac{1}{\sqrt{2}}\frac{1+\epsilon}{1-\epsilon}$. 
So if $R>\frac{1}{\sqrt{2}}$ we can find as much vectors as we want in $B_R$ with mutual distance greater than or equal to $1$ So $R_c({\cal B})\leq \frac{1}{\sqrt{2}}$.\\
So we have $\frac{1}{2}\leq R_c({\cal B})\leq \frac{1}{\sqrt{2}}$ for any Banach space ${\cal B}$.

\end{document}